\documentstyle[amsmath,amssymb,graphicx]{article}

\def\be{\begin{eqnarray}}
\def\ee{\end{eqnarray}}

\def\be{\begin{eqnarray}}
\def\ee{\end{eqnarray}}

\begin{document}

\hfill ITEP/TH-78/07

\centerline{\Large{Gluing of Surfaces with Polygonal Boundaries}}

\bigskip

\centerline{\it E.T.Akhmedov and Sh.Shakirov}

\bigskip

\centerline{ITEP, B.Cheremushkinskaya, 25, Moscow, Russia 117218}

\centerline{Moscow Institute of Physics and Technology,
Dolgoprudny, Russia}

\bigskip

\centerline{ABSTRACT}

\bigskip

By pairwise gluing of edges of a polygon, one produces
two-dimensional surfaces with handles and boundaries. In this
paper, we count the number ${\cal N}_{g,L}(n_1, n_2, \ldots, n_L)$
of different ways to produce a surface of given genus $g$ with $L$
polygonal boundaries with given numbers of edges $n_1, n_2,
\ldots, n_L$. Using combinatorial relations between graphs on real
two-dimensional surfaces, we derive recursive relations between
${\cal N}_{g,L}$. We show that Harer-Zagier numbers appear as a
particular case of ${\cal N}_{g,L}$ and derive a new explicit
expression for them.

\section {Introduction}

A classical question in enumerative combinatorics is: how many
ways are there to glue pairwise all edges of a $2N$-gon so as to
produce a surface of given genus $g$? Such a (\emph{complete})
gluing means that one performs exactly $N$ contractions, so that
no edges of a polygon remain unglued. This problem has been solved
in \cite{HaZ} and the answer is given by the so called
Harer--Zagier numbers $\epsilon_{g}(N)$
\begin{center}
\begin{tabular}{c|c|c|c|cc}
$\epsilon_{g}(N)$ & $g = 0$ & $g = 1$ & $g = 2$ & \ldots \\
\hline
$N = 1$ & 1 & -- & -- & -- \\
\hline
$N = 2$ & 2 & 1 & -- & -- \\
\hline
$N = 3$ & 5 & 10 & -- & -- \\
\hline
$N = 4$ & 14 & 70 & 21 & -- \\
\hline
$N = 5$ & 42 & 420 & 483 & -- \\
\hline
\ldots & \ldots & \ldots & \ldots & \ldots \\
\end{tabular}
\end{center}
with the general formula: \be \epsilon_g (N) =
\frac{(2N)!}{(N+1)!\, (N-2\,g)!}\, \times {\rm Coefficient\,\,
of}\, x^{2g} \,{\rm in}\, \left(\frac{x/2}{\tanh
x/2}\right)^{N+1}.\ee In this paper, we study more general
\emph{incomplete} gluings, i.e, we allow some edges of a polygon
to remain unglued. In this way we generalize the problem and, as a
result, we obtain a simple recursion relation for the
corresponding numbers of the ways to glue.

First of all, let us give a definition of gluing. A particular
gluing of a polygon's boundary is defined via a semi--Gaussian
word along its boundary. A semi--Gaussian word is such a sequence
of letters in which every entry appears \emph{no more} than two
times. This sequence of letters is assigned to the boundary of the
polygon in such a way that one letter corresponds to one edge at
the boundary. If some letter appears twice in the semi--Gaussian
word, then we identify corresponding two edges at the boundary of
the polygon in a way which respects the orientability of the
two--dimensional surface\footnote{In this paper we always consider
orientable two--dimensional surfaces.}. Thus, a complete gluing is
such a word, in which every entry appears exactly two times.

After the gluing --- i.e, after identifying all the edges of a
polygon that correspond to the repeated letters --- we arrive at
the orientable two--dimensional surface, which in general has
polygonal boundaries. At the boundaries we have un--glued edges
corresponding to the letters which appear only one time in the
corresponding semi--Gaussian word. We can have 0--gons (punctured
points), 1--gons, 2--gons, 3--gons (triangles) and etc. among the
boundaries. Note, that if one deals with complete gluings, there
can be only punctured points.

It is worth mentioning now that in this paper we distinguish all
the non--repeated letters in the semi--Gaussian word and, hence,
we distinguish all the edges of the polygonal boundaries of the
resulting two--dimensional surfaces. As well we always consider
such semi--Gaussian words which correspond to the gluings leading
to the surfaces with some fixed in advance sequence of letters on
the edges of the polygonal boundaries.

Now, we consider two gluings of a polygon as equivalent if the
corresponding semi--Gaussian words can be mapped to each other via
a rotation of the polygon and/or via a re--naming of the repeated
letters in the words (see fig.1). Note that such an identification
of the gluings is different from the one adopted in \cite{HaZ}. We
discuss this difference in a greater detail in the third section.

\begin{figure}[t]
\begin{center} \includegraphics[width=230pt]{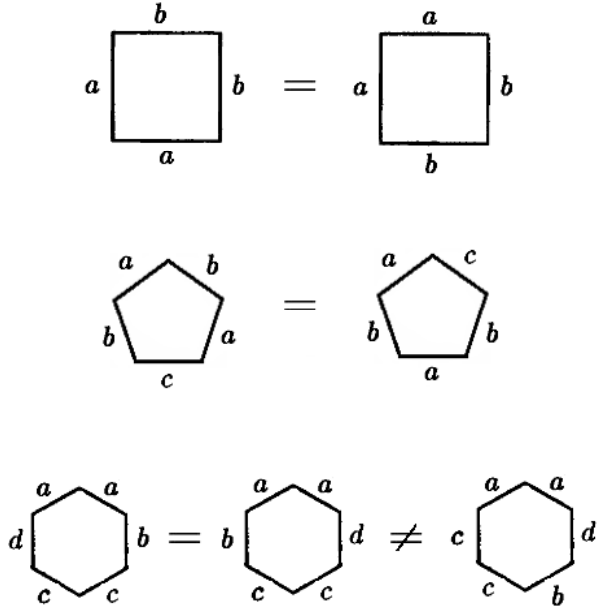}
\caption{Several examples of equivalent and non-equivalent gluings.}
\end{center}
\end{figure}

We use such an equivalence relation for the gluings, because it
corresponds to the natural identification of graphs on
two--dimensional surfaces. Obviously, the image of the polygon's
boundary after the gluing is a graph, embedded into the surface.
The edges of the graph either belong to the polygonal boundaries
or join the vertices belonging to different components of the
boundary. Edges of the graph do not intersect each other, unless
they end at the same vertex of the boundary.

Consider now graphs on two--dimensional orientable surfaces. We
consider two such graphs as equivalent if they can be mapped onto
each other via a diffeomorphism of the surface which respects its
orientation and acts trivially at the boundary (but does not have
to be isotopic to a trivial diffeomorphism). Cutting a surface
along an embedded graph, one will obtain a collection of
two--dimensional surfaces. In this paper we are interested only in
\emph{one--face graphs}: it is such a graph, cutting the surface
along which, one obtains a single polygon without any graphs drawn
on it (except the one going along its boundary). Thus, we arrive
at the statement that there is a one--to--one map between the set
of equivalence classes of all possible gluings of all possible
polygons and the set of equivalence classes of all possible
one--face graphs on all possible two--dimensional surfaces. Such a
one--to--one correspondence is, in fact, a property of the
identification of gluings, adopted in our paper.

 The question that we address is: how many ways are there to glue
edges of an $N$--gon so as to produce a surface of given genus $g$
with $L$ polygonal boundaries with given numbers of edges
$n_1, n_2, \ldots, n_L$? We denote this number as ${\cal N}_{g,L}(n_1, n_2, \ldots, n_L)$.

As one can see, ${\cal N}$ is expressed in terms of the resulting
surface, that is, $g$, $L$ and $n_1, \ldots, n_L$. One naturally
expects ${\cal N}_{g,L}$ to be a symmetric function of $n_1, n_2,
\ldots, n_L$ because the enumeration of holes is arbitrary. Unlike
\cite{HaZ}, we do not treat the total number $N$ of edges of the
glued $N$--gon as an independent parameter: in fact, it can be
easily seen that $N = n_1 + n_2 + \ldots + n_L + 4g + 2L - 2$.
Obviously, in our case $N$ is not necessarily even, namely, we
glue $N$--gon rather than $2N$--gon.

 In this paper we find an explicit formula for ${\cal
N}_{g,L}(n_1, n_2, \ldots, n_L)$. We find ${\cal N}_{g,L}$
inductively, via a recursive relation which expresses the number
of one--face graphs on the surface with $L$ holes and $g$ handles
through the numbers of one--face graphs on the surface with fewer
holes or fewer handles and more holes. This relation has a clear
combinatorial origin following from consecutive cutting of the
surface along the edges of the surface graph.

 As one could expect, ${\cal N}_{g,L}(n_1, n_2, \ldots, n_L)$
reduces to $\epsilon_{g}(N)$ in certain limit. It produces several
explicit formulas for $\epsilon_{g}(N)$, which (to the best of our
knowledge) were not known before.

\section{The calculation of ${\cal N}_{g,L}(n_1, n_2, \ldots, n_L)$}

We count the number of different gluings by means of counting the
number of different graphs. We consider a surface with $g$
handles, which is equipped with $L$ polygonal boundaries, with
numbers of edges $n_1, n_2, \ldots, n_L$. To begin, we assume that
all $n_i > 0$. The case when some $n_i = 0$ will be considered
later.

\paragraph{Lemma I.}
If all $n_i > 0$, then ${\cal N}$ satisfies the following
recursive equation:

\be
\begin{tabular}{c}
$(L + 2g - 1) \cdot {\cal N}_{g,L}(n_1, n_2, \ldots, n_L) = \sum\limits_{i < j} n_i n_j \cdot {\cal N}_{g,L-1}(n_i + n_j + 2, \ \ n_1, \ldots, {\check n}_i, \ldots, {\check n}_j, \ldots, n_L) $ \\
\\
$ + 1/2 \sum\limits_{i} n_i \cdot \sum\limits_{x = 1}^{n_i + 1} {\cal N}_{g-1,L+1}(n_i + 2 - x, x, \ \ n_1, \ldots, {\check n}_i, \ldots, n_L)$,
$1\leq i,j \leq L$. \\
\end{tabular}
\label{MainRecursion} \ee The initial condition for the recursion
is that ${\cal N}_{0,1}(n) = const = 1$, which is obvious, because
a sphere with one hole is already a polygon. If $L \leq 0$ or $g <
0$, it is assumed that ${\cal N} = 0$. The check-sign ${\check n}$
stands for omitting of the arguments.

\paragraph{Proof.} The recursion relation we consider and its derivation
are similar to the ones given in \cite{HaZ} (see page 475). ${\cal
N}_{g,L}$ counts the number of one--face graphs. All one--face
graphs on the same surface consist of the same number of edges. It
can be easily seen that there are $L + 2g - 1$ edges that connect
boundary vertices and do not belong to the boundary. Let us take a
graph from the family ${\cal N}_{g,L}$ and cut the surface along
one of such edges. By the cutting we mean that we unglue two of
the edges of the polygon and assign different letters to them.
I.e., we distinguish two new edges of the resulting surface.

Number ${\cal N}_{g,L}$ of one--face graphs on original surface is
proportional to the sum over the different edges of the numbers of
one--face graphs on the surfaces, obtained by cutting along every
edge. The coefficient of the proportionality is $L + 2g - 1$,
because by doing this, we count one graph on the original surface
exactly $L + 2g - 1$ times. This explains the left hand side of
(\ref{MainRecursion}).

Via the cutting along an edge of a graph from the class ${\cal
N}_{g,L}$ we obtain either graph from the class ${\cal
N}_{g,L-1}$ (i.e. the cutting decreases the number of holes) or
from the class ${\cal N}_{g-1,L+1}$ (i.e. the cutting decreases
the number of handles and increase the number of holes). The
latter two classes correspond to the two terms on the right hand
side of (\ref{MainRecursion}). What is the same, edges of the one--face
graphs are divided into two classes: connecting two different
boundaries, e.g.
\begin{center}
\includegraphics[width=300pt]{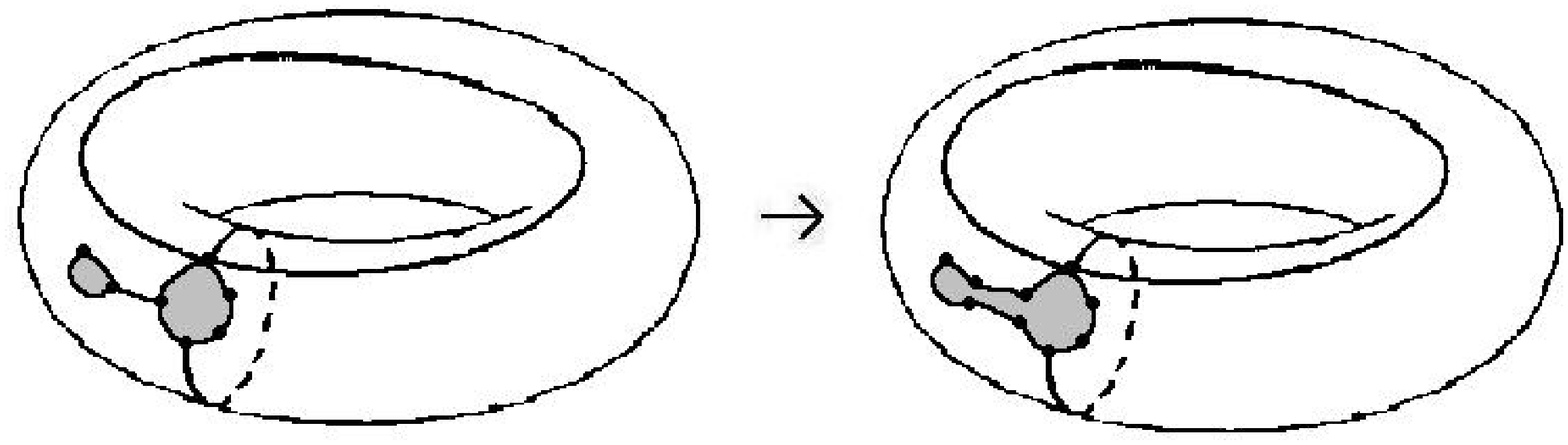}
\end{center}
and connecting one boundary with itself, e.g.
\begin{center}
\includegraphics[width=300pt]{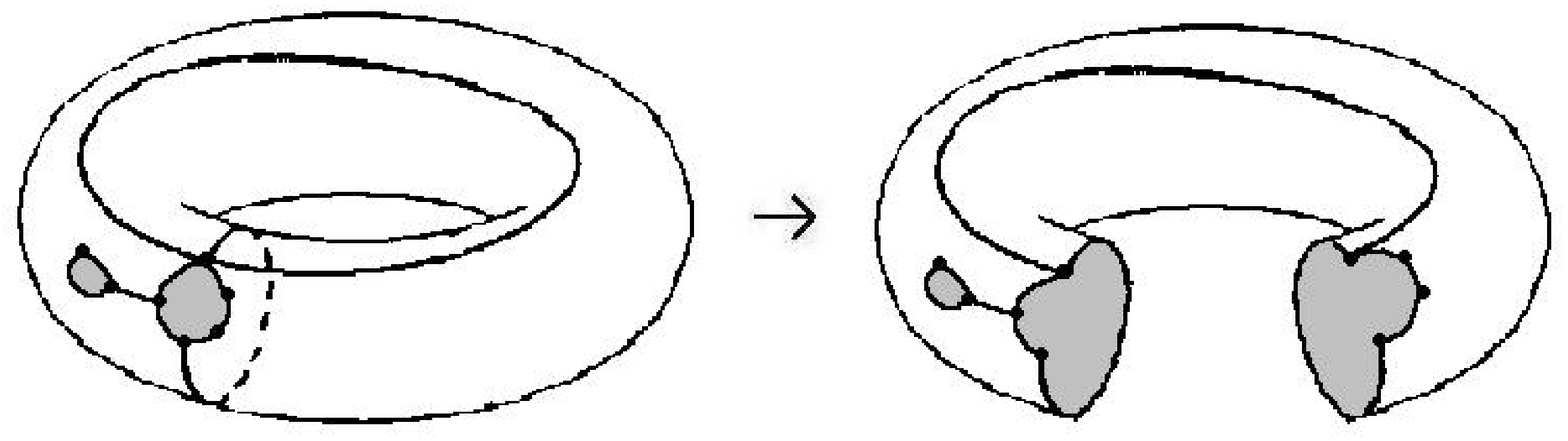}
\end{center}
These two classes correspond to the first and the second terms on
the right hand side of (\ref{MainRecursion}).

\paragraph {$\bullet$} To connect two different (e.g. $i$-th and $j$-th)
boundaries, one can draw an edge from one of $n_i$ vertices of the
first boundary to one of the $n_j$ vertices of the second
boundary. Hence, there are $n_i\,n_j$ non--diffeomorphic to each
other ways to draw one edge connecting the two boundaries in
question. By cutting over such an edge, we obtain a surface with
fewer holes and the same genus. Two holes with $n_i$ and $n_j$
edges transform into a single hole with $n_i + n_j + 2$ edges (see
the first figure). This observation leads to the first term on the
right hand side of (\ref{MainRecursion}).

\paragraph {$\bullet$} To connect an $i$-th boundary with itself, we have to draw
one of the edge encircling one of $g$ handles (otherwise, this
will fail to be a one-face graph). The cut over such an edge
eliminates the corresponding handle and transforms the hole into
two holes (see the second figure). The sum of the numbers of edges
on the new holes is equal to $n_i + 2$, because two new edges
appear after performing the cut.

The boundary's vertices are enumerated from $1$ to $n_i$, in a
natural order. Suppose that the edge terminates in the $p$-th and
$q$-th vertices on the boundary. Such an edge divides the boundary
into two parts. Numbers of edges on these two parts are equal to
$|p - q|$ and $n_i - |p - q|$. After performing the cut, these two
parts transform into two independent holes on the cutted surface.
These holes will have $|p - q| + 1$ and $n_i - |p - q| + 1$ edges.
One finally takes a sum over all edges, that is, over $1 \leq p
\leq q \leq n_i$:
\begin{center}
$ \sum\limits_{p \leq q} {\cal N}_{g-1,L+1}(|p - q| + 1, n_i - |p
- q| + 1, \ \ n_1, \ldots, {\check n}_i, \ldots, n_L) = $

$ 1/2 \sum\limits_{p = 1}^{n_i} \sum\limits_{q = 1}^{n_i} {\cal
N}_{g-1,L+1}(|p - q| + 1, n_i - |p - q| + 1, \ \ n_1, \ldots,
{\check n}_i, \ldots, n_L) + $

$ + 1/2 \sum\limits_{p = q = 1}^{n_i} {\cal N}_{g-1,L+1}(|p - q| +
1, n_i - |p - q| + 1, \ \ n_1, \ldots, {\check n}_i, \ldots, n_L).
$
\end{center}
This expression is simplified by introducing a new summation
variable $x = |p - q| + 1$. Then, it turns into
\begin{center}
$ n_i /2 \sum\limits_{x = 1}^{n_i} {\cal N}_{g-1,L+1}(x, n_i + 2 -
x, \ \ n_1, \ldots, {\check n}_i, \ldots, n_L) + n_i /2 \ {\cal
N}_{g-1,L+1}(1, n_i + 1, \ \ n_1, \ldots, {\check n}_i, \ldots,
n_L) = $

$ n_i / 2 \cdot \sum\limits_{x = 1}^{n_i + 1} {\cal
N}_{g-1,L+1}(x, n_i + 2 - x, \ \ n_1, \ldots, {\check n}_i,
\ldots, n_L).$
\end{center}
This calculation justifies the second term of the right hand side
of (\ref{MainRecursion}).

The lemma is proved.

\paragraph{Theorem I.}
The unique solution of (\ref{MainRecursion}) is given by the
following explicit formula:

\be {\cal N}_{g,L}(n_1, \ldots, n_L) = \dfrac{1}{4^g} \cdot n_1
\ldots n_L \cdot \dfrac{(\Sigma n + 4g + 2L - 3)!}{(\Sigma n + 2g
+ L - 1)!} \cdot \sum\limits_{\lambda_1 + \ldots + \lambda_L = g}
\prod\limits_{k = 1}^{L} \dfrac{(2\lambda_k + n_k)!}{(n_k)!
(2\lambda_k + 1)!} \label{MainFormula} \ee where the sum is taken
over all ordered decompositions $g = \lambda_1 + \ldots +
\lambda_L$ of a non--negative integer number $g$ into $L$
non-negative integers $\lambda_1, \ldots, \lambda_L$. We use the
notation $\Sigma n = n_1 + \ldots + n_L$.

\paragraph{Proof.}

The theorem is proved straightforwardly by induction.
\\
\\
Let us give a proper generalization of (\ref{MainRecursion}) and
(\ref{MainFormula}) to the case when some of $n_i$'s can be zero.
Actually, in this paper, we always assume that there is at least
one unglued edge, i.e. we always consider $\sum_{i=1}^L n_i \geq
1$. Our formulas do not work for the case when all $n_i$ are equal
to zero, i.e. ${\cal N}_{g,L}(0, 0, \ldots, 0)$.

\paragraph{Lemma II.} If some (but not all) $n_i$ are zero, then ${\cal N}$ satisfies
the following recursive equation:

\be
\begin{tabular}{c}
$(L + 2g - 1) \cdot {\tilde {\cal N}}_{g,L}(n_1, n_2, \ldots, n_L) =
\sum\limits_{i < j} {\tilde n}_i {\tilde n}_j \cdot {\tilde {\cal N}}_{g,L-1}(n_i + n_j + 2, \ \ n_1, \ldots, {\check n}_i, \ldots, {\check n}_j, \ldots, n_L) $ \\
\\
$ + 1/2 \sum\limits_{i} {\tilde n}_i \cdot \sum\limits_{x = 1}^{n_i + 1} {\tilde {\cal N}}_{g-1,L+1}(n_i + 2 - x, x, \ \ n_1, \ldots, {\check n}_i, \ldots, n_L)$, \\
\end{tabular}
\label{MainRecursionZeros} \ee where ${\cal N}_{g,L}(n_1, n_2,
\ldots, n_L) = {\tilde {\cal N}}_{g,L}(n_1, n_2, \ldots, n_L) /
\#(n_i = 0)!$ and ${\tilde n} = n$ if $n > 0$ and ${\tilde n} = 1$
if $n = 0$. Here $\#(n_i = 0)$ is the number of $n_i$'s which are
equal to zero.

\paragraph{Proof.} The proof repeats that of the Lemma I, with minor
corrections. A boundary with $0$ edges is a punctured point. It
follows immediately, that one can draw an edge in ${\tilde n}_i
{\tilde n}_j$ ways to connect two different ($i$-th and $j$-th)
boundaries. Finally, it is customary to assume punctured points to
be indistinguishable. Therefore, to avoid counting equivalent
graphs several times, one should divide over the factorial of the
number of punctured points, that is $\#(n_i = 0)!$.

\paragraph{Theorem II.} The unique solution of (\ref{MainRecursionZeros})
is given by the following explicit formula:

\be {\cal N}_{g,L}(n_1, \ldots, n_L) = \dfrac{1}{4^g}
\dfrac{1}{\#(n_i = 0)!} \cdot {\tilde n}_1 \ldots {\tilde n}_L
\cdot \dfrac{(\Sigma n + 4g + 2L - 3)!}{(\Sigma n + 2g + L - 1)!}
\cdot \sum\limits_{\lambda_1 + \ldots + \lambda_L = g}
\prod\limits_{k = 1}^{L} \dfrac{(2\lambda_k + n_k)!}{(n_k)!
(2\lambda_k + 1)!}, \label{MainFormulaZeros} \ee where the sum is
taken over all ordered decompositions $g = \lambda_1 + \ldots +
\lambda_L$ of the integer number $g$ into $L$ non-negative
integers $\lambda_1, \ldots, \lambda_L$.

\paragraph{Proof.} Completely similar to that of the Theorem I.

\subsection{Example: Sphere}
In the $g = 0$ case the surface that is produced is a sphere with
polygonal boundaries. The big decomposition factor in
(\ref{MainFormula}) trivialises, and the formula
(\ref{MainFormula}) reduces to (all $n_i > 0$)

\be {\cal N}_{0,L}(n_1, n_2, \ldots, n_L) = {n}_1 \ldots  {n}_L
\cdot \dfrac{(\Sigma n + 2L - 3)!}{(\Sigma n + L - 1)!}.
\label{SphereFormula} \ee For example,
\begin{center}
$ {\cal N}_{0,1}(n_1) = 1, $ \\
\end{center}

\begin{center}
$ {\cal N}_{0,2}(n_1, n_2) =  {n}_1 \,  {n}_2, $ \\
\end{center}

\begin{center}
${\cal N}_{0,3}(n_1, n_2, n_3) =  {n}_1 \,  {n}_2 \,  {n}_3 \cdot (n_1 + n_2 + n_3 + 3), $ \\
\end{center}

\begin{center}
$ {\cal N}_{0,4}(n_1, n_2, n_3, n_4) =  {n}_1\,  {n}_2 \,  {n}_3
\,  {n}_4 \cdot (n_1 + n_2 + n_3 + n_4 + 4) (n_1 + n_2 + n_3 + n_4
+ 5), $
\end{center}

\begin{center}
$ {\cal N}_{0,5}(n_1, n_2, n_3, n_4, n_5) =  {n}_1 \,
 {n}_2\,  {n}_3 \,  {n}_4\,  {n}_5 \cdot (n_1
+ n_2 + n_3 + n_4 + n_5 + 5) (n_1 + n_2 + n_3 + n_4 + n_5 + 6)
(n_1 + n_2 + n_3 + n_4 + n_5 + 7).$
\end{center}

\subsection{Example: Torus}

In the $g = 1$ case the surface that is produced is a torus with
polygonal boundaries. The big decomposition factor in
(\ref{MainFormula}) now contributes, and the formula
(\ref{MainFormula}) takes the form (all $n_i > 0$)

\be {\cal N}_{1,L}(n_1, n_2, \ldots, n_L) = \dfrac{1}{4}
 {n}_1 \ldots  {n}_L \cdot \dfrac{(\Sigma n + 2L +
1)!}{(\Sigma n + L + 1)!} \cdot \left[ \dfrac{(n_1+1)(n_1+2)}{6} +
\ldots + \dfrac{(n_L+1)(n_L+2)}{6} \right]. \label{TorusFormula}
\ee For example,
\begin{center}
${\cal N}_{1,1}(n_1) =  {n}_1 \, (n_1 + 1) (n_1 + 2) (n_1 + 3)  / 4! $ \\
\end{center}

\begin{center}
${\cal N}_{1,2}(n_1,n_2) =  {n}_1\, {n}_2 \cdot
(n_1 + n_2 + 4) (n_1 + n_2 + 5) (n_1^2 + n_2^2 + 3 n_1 + 3 n_2 + 4) / 4!,$ \\
\end{center}

\begin{center}
${\cal N}_{1,3}(n_1,n_2,n_3) =  {n}_1 \,  {n}_2 \,  {n}_3 \cdot
(n_1 + n_2 + n_3 + 5) (n_1 + n_2 + n_3 + 6) (n_1 + n_2 + n_3 + 7)
(n_1^2 + n_2^2 + n_3^2 + 3 n_1 + 3 n_2 + 3 n_3 + 6) / 4!, $ \\
\end{center}

\begin{center}
${\cal N}_{1,4}(n_1,n_2,n_3,n_4) =  {n}_1\, {n}_2 \,{n}_3 \,{n}_4 \cdot
(n_1 + n_2 + n_3 + n_4 + 6) (n_1 + n_2 + n_3 + n_4 + 7)
(n_1 + n_2 + n_3 + n_4 + 8)(n_1 + n_2 + n_3 + n_4 + 9)
(n_1^2 + n_2^2 + n_3^2 + n_4^2 + 3 n_1 + 3 n_2 + 3 n_3 + 3 n_4 + 8)/4!. $ \\
\end{center}

\section {The Harer-Zagier case}

This section is devoted to comparison with the results of
\cite{HaZ}. We claim, that ${\cal N}_{g,L}(n_1, n_2, \ldots, n_L)$
is a multi-parametric generalization of numbers $\epsilon_{g}(N)$.
Let us specify precisely, what particular case corresponds to
$\epsilon_{g}(N)$.

Gluing of edges in a polygon naturally possess a cyclic symmetry.
For example, the following two complete gluings of a 4-gon into a
sphere
\begin{center}
\includegraphics[width=150pt]{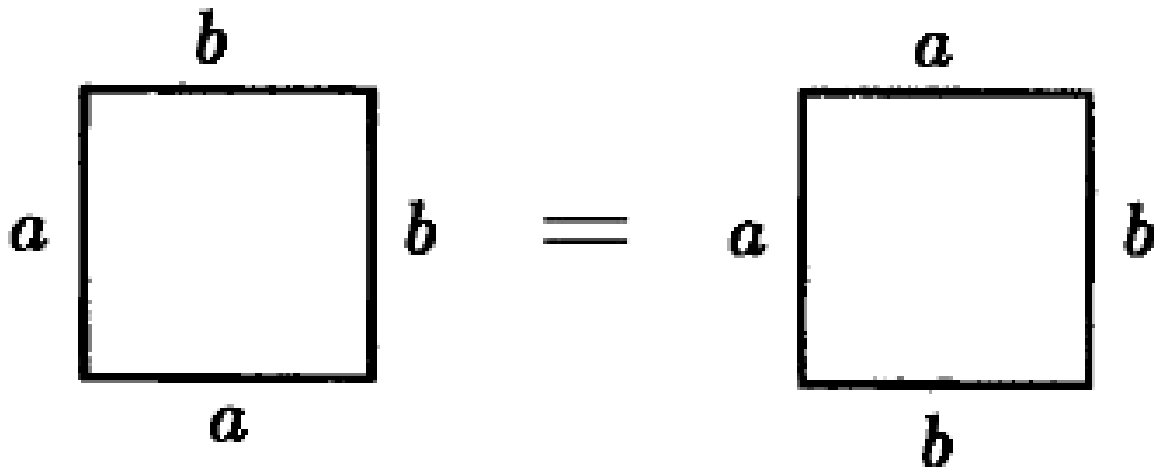}\
\end{center}
are undistinguishable, because they are related to each other by
cyclic exchange of the edges of the square and lead to graphs on
the sphere which are related to each other via a diffeomorphism of
the sphere.

For the problem considered in \cite{HaZ}, it was necessary to deal
with complete gluings of polygons {\emph{without} cyclic symmetry
of edges. The simpliest way to break the cyclic symmetry is to
introduce one distinguished vertex of the polygon:
\begin{center}
\includegraphics[width=150pt]{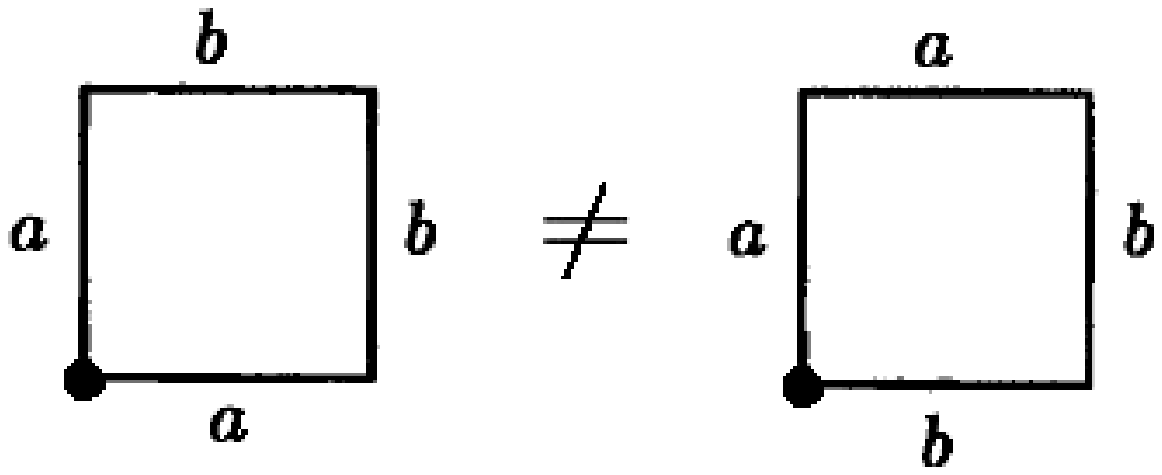}\
\end{center}
Then, two gluings in the above example are different. As a result,
$\epsilon_{0}(2) = 2$. However, this way of breaking the cyclic
symmetry immediately leads to violation of the one-to-one map
between the set of equivalence classes of polygon gluings and the
set of equivalence classes of graphs on the surfaces. For example,
two different gluings
\begin{center}
\ \ \ \ \ \ \ \includegraphics[width=175pt]{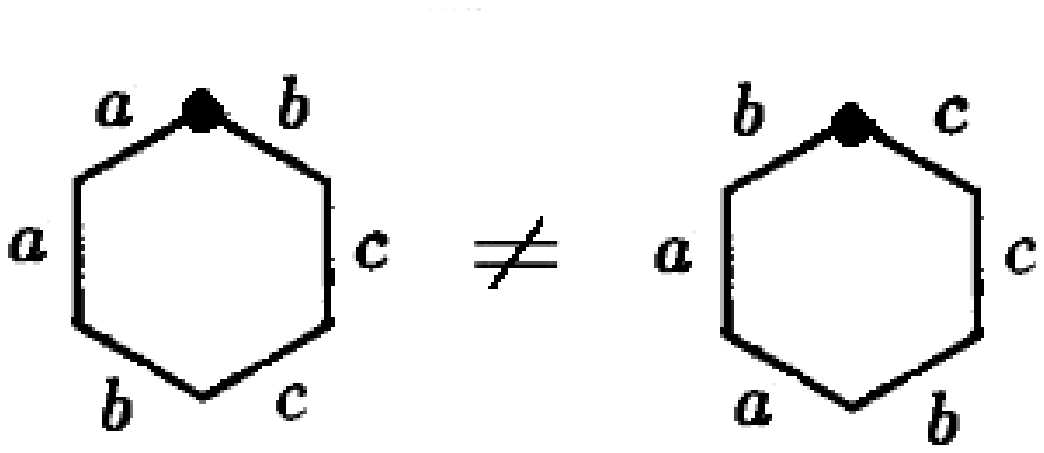}\
\end{center}
lead to the equivalent graphs on the sphere:
\begin{center}
\ \ \ \ \ \ \ \includegraphics[width=175pt]{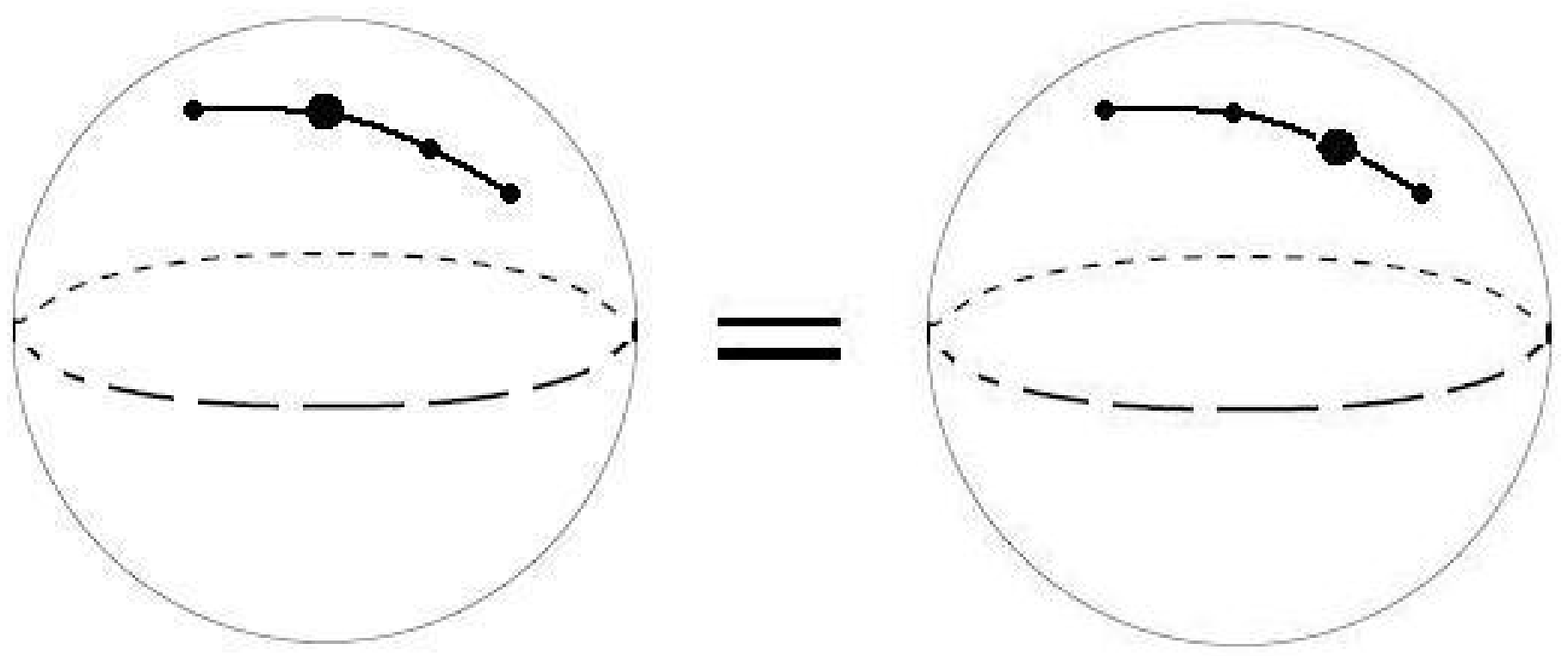}\
\end{center}
In our problem, we would like to respect the one--to--one map in
question. That is why in this paper we do not select any
distinguished point of the polygon. Then, complete gluing numbers
${\cal N}_{g,L}(0, 0, \ldots, 0)$ do not coincide with the
Harer-Zagier numbers.

Instead, it is easy to see that the Harer-Zagier numbers should
coincide with ${\cal N}_{g, L}(1, 0, \ldots, 0)$. The presence of
single unglued edge breaks the symmetry of glued edges under
cyclic exchange. A single unglued edge acts exactly as a
distinguished vertex on the polygon's boundary.

For example, two gluings
\begin{center}
\includegraphics[width=150pt]{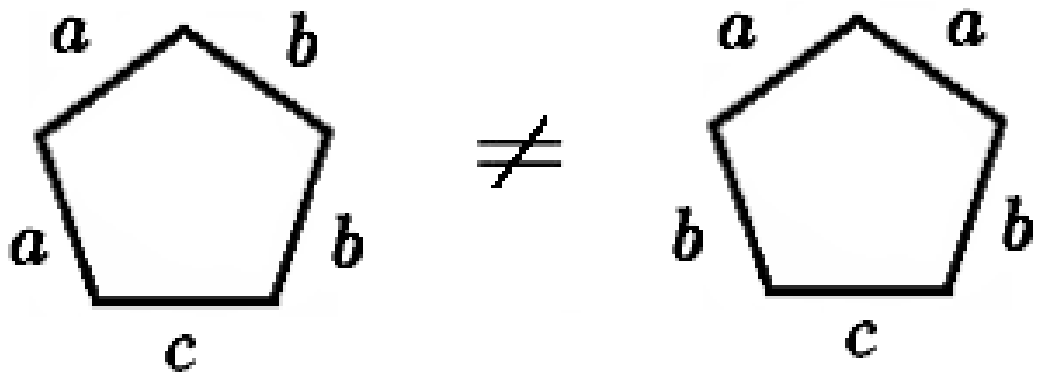}\
\end{center}
are different, therefore ${\cal N}_{0,3}(1, 0, 0) =
\epsilon_{0}(2)$. Generally, the following equality holds:

\be {\cal N}_{g,L}(1, 0, \ldots, 0) = \epsilon_{g}(2g + L - 1)
\label{HZ} \ee In contrast with the distinguished vertex, in the
presence of one unglued edge the one--to--one map in question is
preserved. For example, two different gluings
\begin{center}
\includegraphics[width=175pt]{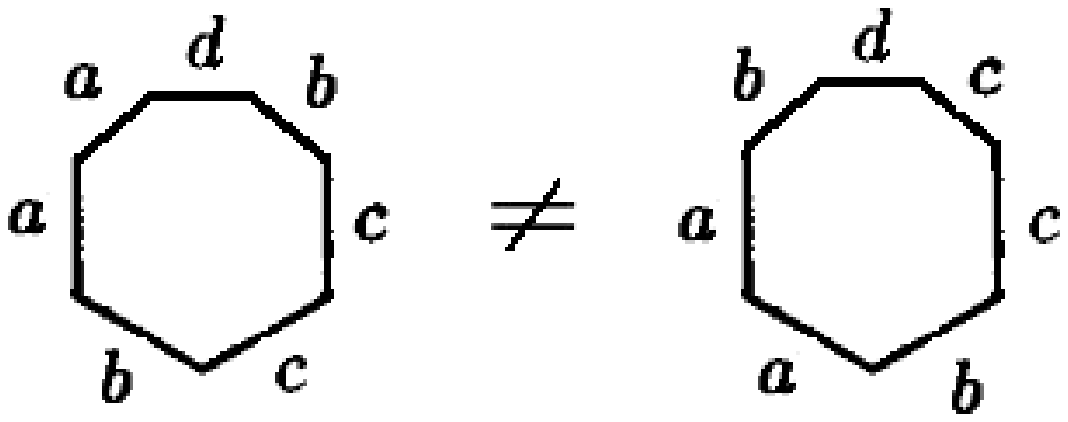}\
\end{center}
lead to two different graphs on the sphere:
\begin{center}
\includegraphics[width=175pt]{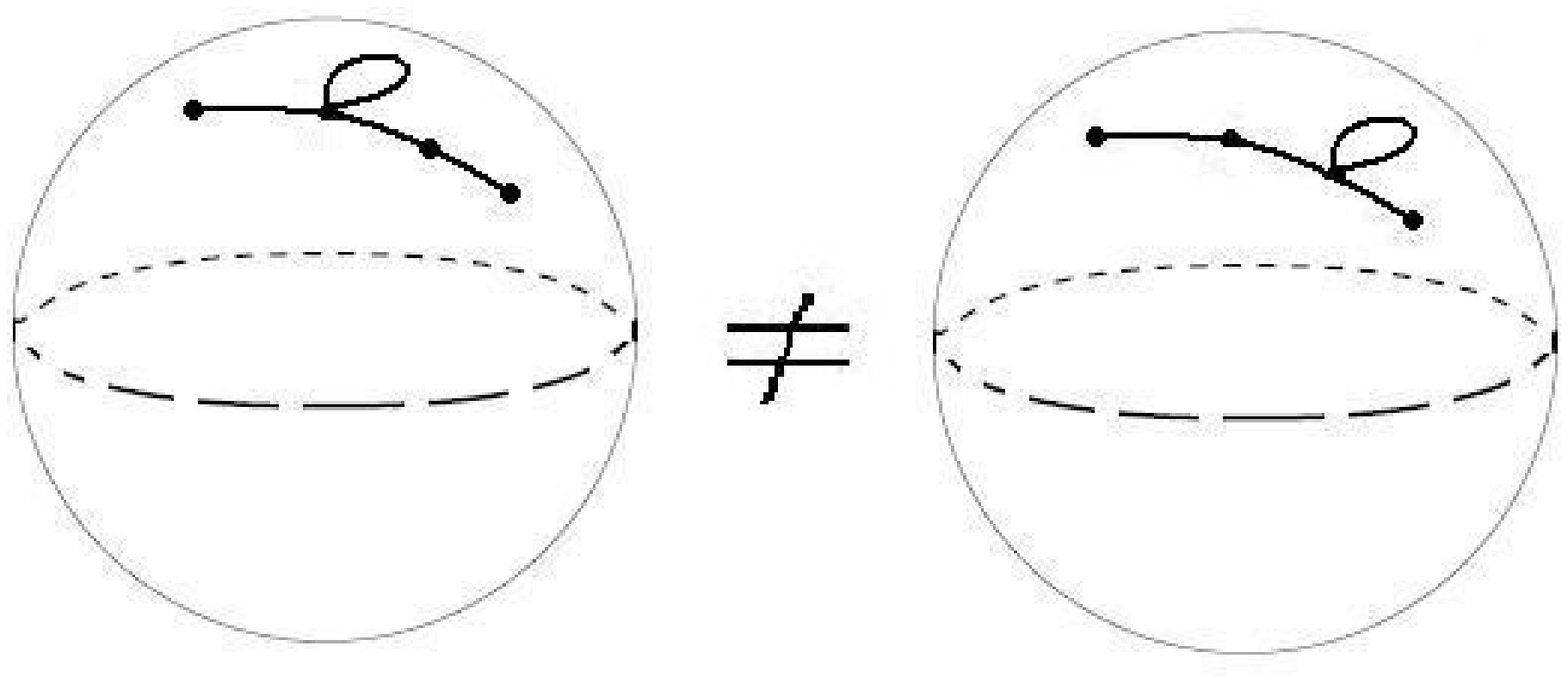}\
\end{center}
Let us make use of the explicit formula for ${\cal N}_{g,L}(n_1,
n_2, \ldots, n_L)$. By substituting (\ref{MainFormulaZeros}) into
(\ref{HZ}) and performing the algebra, one obtains an explicit
representation of the Harer-Zagier's numbers as the following sum
over decompositions:

\be \epsilon_{g}(N) = {\cal N}_{g,N - 2g + 1}(1, 0, \ldots, 0) =
\dfrac{1}{4^g} \cdot \dfrac{(2N)!}{(N - 2g + 1)! (N)!} \cdot
\sum\limits_{\lambda_1 + \ldots + \lambda_L = g}
\dfrac{1}{(2\lambda_1 + 1) \ldots (2\lambda_L + 1)},
 \label{RepHZ}
\ee where $L = N - 2g + 1$. For $g = 0$, the decomposition factor
is trivial (equal to $1$) and we recover the Catalan numbers:

\be \epsilon_{0}(N) = \dfrac{(2N)!}{(N + 1)! (N)!} \cdot
\sum\limits_{\lambda_1 + \ldots + \lambda_L = 0}
\dfrac{1}{(2\lambda_1 + 1) \ldots (2\lambda_L + 1)} =
\dfrac{(2N)!}{(N + 1)! (N)!} = \dfrac{1}{N+1} \big( ^{2N}_{N}
\big).
 \label{CatalanHZ}
\ee
For $g = 1$, the sum can be easily calculated, because all
decompositions of $1$ into $L$ parts have a form
$[0,\ldots,1,\ldots,0]$. In this case, the decomposition factor is
equal to $L/3$ and we obtain

\be \epsilon_{1}(N) = \dfrac{(2N)!}{4(N - 1)! (N)!} \cdot
\sum\limits_{\lambda_1 + \ldots + \lambda_L = 1}
\dfrac{1}{(2\lambda_1 + 1) \ldots (2\lambda_L + 1)} = \dfrac{(2N)!
\cdot (N-1)/3}{4(N - 1)! (N)!} = \dfrac{1}{12}
\dfrac{(2N)!}{(N-2)!N!}.
 \label{ToricHZ}
\ee Similar calculations can be done for $g > 1$. The well-known
Harer-Zagier generating function

$$1 + 2 \sum\limits_{g = 0}^{\infty} \sum\limits_{N = 2g}^{\infty} \epsilon_{g}(N) x^{N+1} y^{N - 2g + 1}
\dfrac{1}{(2N - 1)!!} = \left(\dfrac{1 + x}{1 - x}\right)^y$$
follows from from our formula (\ref{RepHZ}) at once:
\[
\begin{array}{cc}
1 + 2 \sum\limits_{g = 0}^{\infty} \sum\limits_{N = 2g}^{\infty} \epsilon_{g}(N) x^{N+1} y^{N - 2g + 1} \dfrac{1}{(2N - 1)!!} = \\
\\
= 1 + 2 \sum\limits_{g = 0}^{\infty} \sum\limits_{N = 2g}^{\infty} x^{N+1} y^{N - 2g + 1} \dfrac{1}{4^g} \cdot \dfrac{(2N)!}{(2N - 1)!!(N)!} \dfrac{1}{(N - 2g + 1)!} \cdot \sum\limits_{\lambda_1 + \ldots + \lambda_L = g} \dfrac{1}{(2\lambda_1 + 1) \ldots (2\lambda_L + 1)} = \\
\\
= 1 + 2 \sum\limits_{g = 0}^{\infty} \sum\limits_{N = 2g}^{\infty} x^{N+1} y^{N - 2g + 1} 2^{-2g} 2^N \cdot  \dfrac{1}{(N - 2g + 1)!} \cdot \sum\limits_{\lambda_1 + \ldots + \lambda_L = g} \dfrac{1}{(2\lambda_1 + 1) \ldots (2\lambda_L + 1)} = \\
\\
\end{array}
\]
\[
\begin{array}{cc}
= 1 + \sum\limits_{g = 0}^{\infty} \sum\limits_{L = 1}^{\infty} (2x)^{L} y^{L} x^{2g} \cdot  \dfrac{1}{L!} \cdot \sum\limits_{\lambda_1 + \ldots + \lambda_L = g} \dfrac{1}{(2\lambda_1 + 1) \ldots (2\lambda_L + 1)} = \\
\\
= 1 + \sum\limits_{L = 1}^{\infty} \dfrac{(2y)^{L}}{L!} \sum\limits_{g = 0}^{\infty} \sum\limits_{\lambda_1 + \ldots + \lambda_L = g} \dfrac{x^{2g + L}}{(2\lambda_1 + 1) \ldots (2\lambda_L + 1)} = \\
\\
= 1 + \sum\limits_{L = 1}^{\infty} \dfrac{(2y)^{L}}{L!} \left[\sum\limits_{\lambda = 0}^{\infty}
\dfrac{x^{2 \lambda + 1}}{(2\lambda + 1)} \right]^L = \\
\\
= 1 + \sum\limits_{L = 1}^{\infty} \dfrac{(2y)^{L}}{L!}
\left[\dfrac{1}{2} \, \log \dfrac{1 + x}{1 - x} \right]^L =
\exp\left[ 2 y \cdot \dfrac{1}{2} \log \dfrac{1 + x}{1 - x}
\right] = \left(\dfrac{1 + x}{1 - x}\right)^y.
\end{array}
\]
\section{Conclusion}
Instead of conclusion let us explain what construction is behind
the ${\cal N}_{g,L}(n_1, n_2, \ldots, n_L)$ numbers. It is well
known that cellular decompositions of polygons (with fixed
orientation and one distinguished edge) by diagonals are in
one--to--one correspondence with Stasheff polytops (see e.g.
\cite{Stasheff}). Let us explain the correspondence.

One can draw at the same time $N-3$ non--intersecting diagonals in
an $N$--gon. This way we obtain a triangulation of the polygon.
The correspondence in question is pure combinatorial: a
$k$--dimensional ($k\leq N-3$) face of a polytop corresponds to
the cellular decomposition of the polygon by its diagonals, which
can be obtained from the polygon's triangulation via removal of
$k$ among $N-3$ diagonals. A $k$--dimensional face of the Stasheff
polytop belongs to the boundary of an $n$--dimensional face
($N-3\geq n>k$) if and only if the cellular decomposition with
$N-3-n$ diagonals (corresponding to the $n$--dimensional face) can
be embedded into the cellular decomposition with $N-3-k$ diagonals
(corresponding to the $k$--dimensional face). E.g. vertices of the
Stasheff polytop are diagonal triangulations of the polygon, edges
--- are triangulations without one diagonal and etc.

If we recall that the polygon topologically is just a disc, one
would like to consider a generalization of that construction to
real two--dimensional surfaces with handles and polygonal
boundaries. I.e. we would like to build combinatorial complexes
related to cellular decompositions of real two--dimensional
surfaces. So that a $k$--dimensional face of such a complex is
related to the triangulation of the surface with $k$ edges of the
triangulation removed. Obviously these $k$--dimensional faces are
just $k$--dimensional Stasheff polytops. From this construction it
is clear that any one--face graph on the surface is related to the
maximal dimension cell in the corresponding complex. Then, ${\cal
N}_{g,L}(n_1, n_2, \ldots, n_L)$ counts just the number of cells
of the maximal dimension in the complex corresponding to the
surface with $g$ handles and $L$ polygonal holes with $n_1, n_2,
\ldots, n_L$ edges.

 In the literature one usually considers ribbon graphs, which are dual to
the cellular decompositions, on the real two--dimensional surfaces
instead of the cellular decompositions themselves. In that case
the complex in question is referred to as graph complex (see e.g.
\cite{graphcompl}). The important problem is to calculate the
homology of this complex.

There are other interesting questions related to the ${\cal
N}_{g,L}(n_1, n_2, \ldots, n_L)$ numbers: such as finding the
generating function for these numbers and/or matrix integral
representation. Graphs on surfaces appear in various branches of
modern physics, such as matrix models of two-dimensional quantum
gravity \cite{qgrav}, so an existence of a matrix integral
representation for ${\cal N}_{g,L}(n_1, n_2,\ldots, n_L)$ would
not be surprising.

\section{Acknowledgements}

ETA would like to thank I.Artamkin and S.Lando for valuable
discussions and A.Babichev for lessons of mathematical culture and
grammar. As well we would like to thank the referee of
``Functional Analysis and Applications'' journal for very valuable
comments and corrections. This work was partially supported by the
Federal Nuclear Energy Agency. The work of S.S. was partially
supported by the Program for Supporting Leading Scientific Schools
(Grant No. NSh-8004.2006. 2), the Russian Foundation for Basic
Research (Grant Nos. RFBR–Italy 06-01-92059-CE and 07-02-00642)
and by the Dynasty Foundation.

\end{document}